\documentclass[12pt]{article}

\usepackage{amsmath,amsthm,amsfonts,latexsym,amsopn,verbatim,amscd,amssymb}

\theoremstyle{plain}
\newtheorem{theorem}{Theorem}[section]

\newtheorem{proposition}[theorem]{Proposition}
\newtheorem{corollary}[theorem]{Corollary}
\theoremstyle{definition}
\numberwithin{equation}{section}

\DeclareMathOperator{\spec}{Spec}
\DeclareMathOperator{\cent}{Cent}

\newcommand{\bnum}{\begin{enumerate}}
\newcommand{\enum}{\end{enumerate}}

\begin{document}

\title{Spectrum and genus of commuting  graphs of some classes of finite rings }
\author{Jutirekha Dutta  and Rajat Kanti Nath\footnote{Corresponding author}  }
\date{}
\maketitle
\begin{center}\small{ Department of Mathematical Sciences,\\ Tezpur
University,  Napaam-784028, Sonitpur, Assam, India.\\
Emails: jutirekhadutta@yahoo.com, rajatkantinath@yahoo.com}
\end{center}


\smallskip

\noindent {\small{\textbf{Abstract:} }
The commuting graph of a non-commutative ring $R$ with center $Z(R)$ is a simple undirected graph whose vertex set is $R\setminus Z(R)$ and two vertices $x, y$ are adjacent if and only if $xy = yx$.
In this paper, we compute the spectrum and genus of  commuting graphs of some classes of  finite rings. 
 }

\bigskip



\section{Introduction} \label{S:intro}
Let $R$ be a non-commutative ring with center $Z(R)$. The commuting graph of $R$, denoted by $\Gamma_R$, is a simple undirected graph whose vertex set is $R\setminus Z(R)$ and two vertices $x, y$ are adjacent if and only if $xy = yx$. In recent years, many mathematicians have considered  commuting graph of different rings and studied  various graph theoretic aspects (see \cite{a08,aghm04,ar06,ekn15,m10,ov11,vrb14}). Some generalizations of $\Gamma_R$ are also considered in \cite{aBHK15, dBN16}. 
In this paper, we compute the spectrum and genus of the commuting graph of some classes of  finite rings.

Recall that the spectrum of a graph  ${\mathcal{G}}$ denoted by $\spec({\mathcal{G}})$ is the set $\{\lambda_1^{k_1}, \lambda_2^{k_2},$ $\dots, \lambda_n^{k_n}\}$, where $\lambda_1,  \lambda_2, \dots, \lambda_n$ are the eigenvalues of the adjacency matrix of $\mathcal{G}$ with multiplicities $k_1, k_2, \dots, k_n$ respectively. A graph ${\mathcal{G}}$ is called integral if $\spec({\mathcal{G}})$ contains only integers. It is well-known that the complete graph $K_n$ on $n$ vertices is integral and $\spec(K_n) = \{(-1)^{n - 1}, (n - 1)^1\}$. Further, if $\mathcal{G} = \overset{l}{\underset{i = 1}{\sqcup}}K_{m_i}$, where $K_{m_i}$ are complete graphs on $m_i$ vertices for $1 \leq i \leq l$, then 
\begin{equation}\label{spectrum}
\spec(\mathcal{G}) = \{(-1)^{\underset{i = 1}{\overset{l}{\sum}}m_i - l},\, (m_1 - 1)^1,\, (m_2 - 1)^1,\, \dots,\, (m_l - 1)^1\}.
\end{equation} 
The notion of integral graph was introduced by 
Harary and  Schwenk \cite{hS74} in the year 1974. Since then many mathematicians have considered integral graphs, see for example \cite{aV09,iV07,wLH05}. In \cite{dn15,Nath15}, the authors have determined several  groups whose commuting graphs are integral.

The genus of a graph is the smallest non-negative integer $n$ such that the graph can be embedded on the surface obtained by attaching $n$ handles to a sphere. We write $\gamma({\mathcal{G}})$ to denote the genus of a graph ${\mathcal{G}}$. It is worth mentioning that 
\[
\gamma(K_n) = \left\lceil\frac{(n - 3)(n - 4)}{12}\right\rceil \text{ if } n \geq 3 \text{ (see \cite[Theorem 6-38]{wAT73})}.
\] 
Also, if ${\mathcal{G}} = \underset{B\in {\mathcal{B}}}{\sqcup}B$, where ${\mathcal{B}}$ denotes the block decomposition of ${\mathcal{G}}$, then by    \cite[Corollary 1]{bhky62} we have
\begin{equation}\label{genus}
\gamma({\mathcal{G}}) = \underset{B\in {\mathcal{B}}}{\sum}\gamma(B).
\end{equation}

A graph ${\mathcal{G}}$ is called planar or toroidal if $\gamma({\mathcal{G}}) = 0$ or $1$ respectively. As a consequence of our results we show that the commuting graphs of non-commutative rings of order $p^2$ and $p^3$ are integral but not toroidal.

 For any element $r$ of a ring $R$, the set $C_R(r) = \{s \in R : rs = sr\}$ is called the centralizer of  $r $ in $R$. Let $|\cent(R)| = |\{C_R(r) : r \in R\}|$, that is the number of distinct centralizers in $R$. A ring $R$ is called an $n$-centralizer ring if $|\cent(R)| = n$. In \cite{dbn15}, Dutta et al. have  characterized finite $n$-centralizer rings for $n = 4, 5$. As a   consequence of our results,  we  show  that the commuting graphs of   $4, 5$-centralizer finite rings are  integral but not toroidal. Further, we show that the commuting graph of a finite $(p + 2)$-centralizer $p$-ring is  integral for any prime $p$. We conclude this paper by computing the spectrum and genus of the commuting graphs of some finite rings with some specific commuting probability.

\section{ Main results}
A non-commutative ring $R$ is called a CC-ring if all the centralizers of its non-central elements are commutative. In \cite{ekn15}, Erfanian et al. have initiated the study of  CC-rings. In particular, they have computed the diameter of $\Gamma_R^c$ and showed that the clique number and chromatic number of $\Gamma_R^c$ are same for a CC-ring $R$, where $\Gamma_R^c$ denotes the complement of $\Gamma_R$. In the following theorem we compute the spectrum and genus of $\Gamma_R$ for a finite CC-ring $R$.
\begin{theorem}\label{cc-ring}
Let $R$ be a finite CC-ring with distinct centralizers $S_1, S_2, \dots, S_n$  of non-central elements of $R$. Then
\[
\spec(\Gamma_R) = \{(-1)^{\overset{n}{\underset{i = 1}{\sum}}|S_i| - n(|Z(R)| + 1) }, (|S_1| - |Z(R)| - 1)^1,\dots, (|S_n| - |Z(R)| - 1)^1\}
\]
and $\gamma(\Gamma_R) = \overset{n}{\underset{i = 1}{\sum}} \gamma(K_{|S_i| - |Z(R)|})$.
\end{theorem} 
\begin{proof}
Let $R$ be a finite CC-ring and  $S_1, S_2, \dots, S_n$  be the distinct centralizers of  non-central elements of $R$. Let $S_i = C_R(s_i)$ where $s_i \in R\setminus Z(R)$ for $i = 1, 2, \dots, n$. 
Let $s \in  (S_i\cap S_j)\setminus Z(R)$ for some $i, j$ such that $1\leq i \ne j \leq n$. Then $s$ commutes with $s_i$ as well as $s_j$. Let $t \in C_R(s)$ then $ts_i = s_it$ since $s_i \in C_R(s)$ and $R$ is a CC-ring. Therefore $t \in C_R(s_i)$ and so $C_R(s) \subseteq C_R(s_i)$. 
Again, let $u \in C_R(s_i)$ then $us =su$ since $s \in S_i = C_R(s_i)$ and $R$ is a CC-ring. Therefore, $u \in C_R(s)$ and so $C_R(s_i) \subseteq C_R(s)$. Thus $C_R(s) = C_R(s_i)$. Similarly, it can be seen that $C_R(s) = C_R(s_j)$. Hence $C_R(s) = C_R(s_i) = C_R(s_j)$ which is a contradiction. Therefore, $S_i\cap S_j = Z(R)$ for $1\leq i \ne j \leq n$. This shows that $\Gamma_R = \overset{n}{\underset{i = 1}{\sqcup}} K_{|S_i| - |Z(R)|}$. Now the results follow from \eqref{spectrum} and \eqref{genus}.
\end{proof} 

\begin{corollary}\label{CC-cor}
Let $R$ be a finite  CC-ring and $A$ be any finite commutative ring.  Then    the spectrum of the commuting graph of   $R\times A$ is given by
\begin{align*}
 \{ (-1)^{\overset{n}{\underset{i = 1}{\sum}}|A|(|S_i| - |Z(R)|) - n}, (|A|(|S_1| - |Z(R)|) - &1))^1,\dots,\\
&  (|A|(|S_n| - |Z(R)|) - 1))^1 \} 
\end{align*}
and $\gamma(\Gamma_R) = \overset{n}{\underset{i = 1}{\sum}} \gamma(K_{|A|(|S_i| - |Z(R)|)})$
where $S_1,\dots, S_n$ are the distinct centralizers of non-central elements of $R$.  
\end{corollary}

\begin{proof}
Note that $Z(R\times A) = Z(R)\times A$ and $S_1\times A,  S_2 \times A,\dots, S_n \times A$  are the distinct centralizers of non-central elements of $R\times A$.  Therefore, if $R$ is an  CC-ring then $R\times A$ is also a  CC-ring. Hence, the result follows from Theorem \ref{cc-ring}. 
\end{proof}

In general it is difficult to determine all the finite non-commutative rings whose commuting graphs are integral. However, by Theorem \ref{cc-ring}, it follows that the commuting graph of     a finite CC-ring is integral. Further, if $R$ is a finite  CC-ring and $A$ is any finite commutative ring then, by Corollary \ref{CC-cor}, the commuting graph of $R \times A$ is also integral.   In the next result we consider a particular class of CC-rings and compute the spectrum and genus of its commuting graph.
 
\begin{theorem}\label{main2}
Let $R$ be a finite ring such that the additive quotient group $\frac{R}{Z(R)}$ is isomorphic to ${\mathbb{Z}}_p \times {\mathbb{Z}}_p$, where $p$ is a prime. Then 
\[
\spec(\Gamma_R) = \{(-1)^{(p^2 - 1)|Z(R)| - p - 1}, ((p - 1)|Z(R)| - 1)^{p + 1}\} \text{ and } 
\]
\[
\gamma(\Gamma_R) =  (p + 1) \gamma(K_{(p - 1)|Z(R)|}).
\]
\end{theorem}
\begin{proof}
Let  $\frac{R}{Z(R)} = \langle a+Z(R), b+Z(R) : pa, pb \in Z(R)\rangle$, where $a, b \in R$ with $ab \ne ba$. Then for any $z \in Z(R)$, we have
\begin{align*}
C_R(a) &= C_R(ia+z) \\
 &= Z(R) \sqcup a+Z(R) \sqcup \cdots \sqcup (p -1)a+Z(R) \text{ for } 1 \leq i \leq p - 1;
\end{align*}
and   for  $1 \leq j \leq p$,
\begin{align*}
C_R(ja+b) &= C_R(ja+b+z)\\
& = Z(R) \sqcup (ja+b)+Z(R) \sqcup \cdots \sqcup ((p -1)ja+(p -1)b)+Z(R).
\end{align*}
These are the only  centralizers of non-central elements of $R$. Also note that these centralizers are commutative subrings of $R$. Thus $R$ is a CC-ring and
\[
\Gamma_R = K_{|C_R(a)\setminus Z(R)|} \sqcup (\underset{j = 1}{\overset{p}{\sqcup}} K_{|C_R(ja+b)\setminus Z(R)|}).
\]
 That is, $\Gamma_R = K_{(p - 1)|Z(R)|} \sqcup (\underset{j = 1}{\overset{p}{\sqcup}} K_{(p - 1)|Z(R)|}) = \Gamma_R =  \underset{j = 1}{\overset{p + 1}{\sqcup}} K_{(p - 1)|Z(R)}$,  since  $|C_R(a)| = p|Z(R)|$ and $|C_R(ja+b)| = p|Z(R)|$ for $1 \leq j \leq p$.
 Hence the results follow from \eqref{spectrum} and \eqref{genus}.
\end{proof}

\section{Some consequences}
In this section, we obtain several consequences of the results obtained in Section $2$. In general it is difficult to determine all finite non-commutative rings whose commuting graphs are planar or toridal. In this section, we characterize some finite rings whose commuting graphs are planar.  We begin with the following result.
\begin{proposition}\label{App-1}
Let $R$ be a finite ring such that the additive quotient group $\frac{R}{Z(R)}$ is isomorphic to ${\mathbb{Z}}_p \times {\mathbb{Z}}_p$, where $p$ is a prime. Then
\begin{enumerate}
\item $\Gamma_R$ is integral but not toroidal.
\item $\Gamma_R$ is planar if and only if $p = 2, |Z(R)| = 1, 2, 3$ or $4$; or $p = 3, |Z(R)| = 1$ or $2$.
\end{enumerate}
\end{proposition}
\begin{proof}
Part (a) follows from Theorem \ref{main2}. If $p = 2$ then $\gamma(\Gamma_R) = 0$ if and only if $3\gamma(K_{|Z(R)|}) = 0$ if and only if $|Z(R)| = 1, 2, 3$ or $4$. If $p = 3$ then $\gamma(\Gamma_R) = 0$ if and only if $4\gamma(K_{2|Z(R)|}) = 0$ if and only if $|Z(R)| = 1$ or $2$. Hence, part (b) follows.
\end{proof}

\begin{proposition}
Let  $p$ be a prime and $R$ a non-commutative ring of order $p^2$. Then
\begin{enumerate}
\item $\Gamma_R$ is integral but not toroidal.
\item $\Gamma_R$ is planar if and only if $p = 2, 3$ or $5$.
\end{enumerate}
\end{proposition}

\begin{proof}
Note that $|Z(R)| = 1$ and  the additive quotient group $\frac{R}{Z(R)}$ is isomorphic to ${\mathbb{Z}}_p \times {\mathbb{Z}}_p$. So, by Theorem \ref{main2}, we have
\[
\spec(\Gamma_R) = \{(-1)^{p^2 - p - 2}, (p - 2)^{p + 1}\} \text{ and } \gamma(\Gamma_R) =  (p + 1) \gamma(K_{p - 1}).
\]
Thus it follows that $\Gamma_R$ is integral. Also, $\gamma(\Gamma_R) \ne 1$, that is, $\Gamma_R$ is not toroidal. Part (b) follows from the fact that $\gamma(\Gamma_R) = 0$ if and only if $p = 2, 3$ or $5$.
\end{proof}



\begin{proposition}
Let  $p$ be a prime and $R$ a non-commutative ring of order $p^3$. Then
\begin{enumerate}
\item $\Gamma_R$ is integral but not toroidal.
\item $\Gamma_R$ is  planar if and only if $p = 2$.
\end{enumerate}
\end{proposition}

\begin{proof}
Note that $|Z(R)| = p$ and  the additive quotient group $\frac{R}{Z(R)}$ is isomorphic to ${\mathbb{Z}}_p \times {\mathbb{Z}}_p$. So, by Theorem \ref{main2}, we have
\[
\spec(\Gamma_R) = \{(-1)^{p^3 - 2p - 1}, (p^2 - p - 1)^{p + 1}\} \text{ and } \gamma(\Gamma_R) =  (p + 1) \gamma(K_{p^2 - p}).
\]
Thus $\Gamma_R$ is integral. Also, $\gamma(\Gamma_R) \ne 1$, that is, $\Gamma_R$ is not toroidal. Part (b) follows from the fact that $\gamma(\Gamma_R) = 0$ if and only if $p = 2$.
\end{proof}

\begin{proposition}\label{4-cent}
If $R$ is a finite $4$-centralizer ring then $\Gamma_R$ is  integral but not toroidal. Also $\Gamma_R$ is planar if and only if $|Z(R)| = 1, 2, 3$ or $4$.      
\end{proposition}
\begin{proof}
If $R$ is a finite $4$-centralizer ring then by   \cite[Theorem 3.2]{dbn15} we have that the additive quotient group $\frac{R}{Z(R)}$ is isomorphic to ${\mathbb{Z}}_2 \times {\mathbb{Z}}_2$. Now the results follow from Proposition \ref{App-1}.
\end{proof}

\begin{proposition}\label{5-cent}
If $R$ is a finite $5$-centralizer  ring then $\Gamma_R$ is  integral but not toroidal. Also $\Gamma_R$ is planar if and only if $|Z(R)| = 1$ or $2$.
\end{proposition}
\begin{proof}
If $R$ is a finite $5$-centralizer ring then by  \cite[Theorem 4.3]{dbn15} we have that the additive quotient group $\frac{R}{Z(R)}$ is isomorphic to ${\mathbb{Z}}_3 \times {\mathbb{Z}}_3$. Hence, the  result follows from Proposition \ref{App-1}. 
\end{proof}

\noindent In the following proposition, we   compute the spectrum and genus of a $(p+2)$-centralizer $p$-ring, for any prime $p$.

\begin{proposition}
If $R$ is a finite $(p+2)$-centralizer $p$-ring, for any prime $p$, then
\[
\spec(\Gamma_R) = \{(-1)^{(p^2 - 1)|Z(R)| - p - 1}, ((p - 1)|Z(R)| - 1)^{p + 1}\}
\]
and $\gamma(\Gamma_R) =  (p + 1) \gamma(K_{(p - 1)|Z(R)|})$.
\end{proposition}
\begin{proof}
If $R$ is a finite $(p + 2)$-centralizer $p$-ring then by  \cite[Theorem 2.12]{dbn15} the additive quotient group  $\frac{R}{Z(R)}$ is isomorphic to ${\mathbb{Z}}_p \times {\mathbb{Z}}_p$. Hence, the  result follows from Theorem \ref{main2}.
\end{proof}

The commuting probability of a ring $R$ denoted by $\Pr(R)$ is the probability that a randomly chosen pair of elements of $R$ commute. If $R$ is a finite ring then $\Pr(R)$ is given by the ratio
\[
\frac{|\{(r, s) \in R\times R : rs = sr\}|}{|R|^2}.
\]
 The study of $\Pr(R)$ was initiated by MacHale \cite{dmachale} in 1976. MacHale \cite{dmachale} proved the following result.
\begin{theorem}\label{pr-1}
Let $R$ be a finite ring and $p$ the smallest prime divisor of $|R|$. Then $\Pr(R) \leq \frac{p^2 + p - 1}{p^3}$. The equality holds if and only if the additive quotient group $\frac{R}{Z(R)}$ is isomorphic to ${\mathbb{Z}}_p \times {\mathbb{Z}}_p$. 
\end{theorem}

As a consequence of Theorem \ref{main2} and Theorem \ref{pr-1} we have the following results. 
\begin{proposition}
Let $R$ be a finite ring with $\Pr(R) = \frac{5}{8}$ then 
\[
\spec(\Gamma_R) = \{(-1)^{3(|Z(R)| - 1)}, (|Z(R)| - 1)^3\} \text{ and } \gamma(\Gamma_R) =  3\gamma(K_{|Z(R)|}).
\]
\end{proposition}
\begin{proposition}
Let $R$ be a finite ring and $p$ the smallest prime divisor of $|R|$. If $\Pr(R) = \frac{p^2 + p - 1}{p^3}$ then $\spec(\Gamma_R) = \{(-1)^{(p^2 - 1)|Z(R)| - p - 1}, ((p - 1)|Z(R)| - 1)^{p + 1}\}$ and $\gamma(\Gamma_R) =  (p + 1) \gamma(K_{(p - 1)|Z(R)|})$. 
\end{proposition}


\end{document}